\documentclass[11pt]{article}
\usepackage{amsmath,amssymb}

\def\N{\mathbb{N}}
\def\Z{\mathbb{Z}}
\def\R{\mathbb{R}}

\def\C{C^{\infty}(M)}

\newtheorem{definition}{Definition}[section]
\newtheorem{lemma}[definition]{Lemma}
\newtheorem{proposition}[definition]{Proposition}
\newtheorem{theorem}[definition]{Theorem}
\newtheorem{remark}[definition]{Remark}

\newtheorem{example}[definition]{Example}
\newenvironment{proof}{\noindent{\bf Proof.}}{\hfill $\blacklozenge$}

\begin{document}

\title{On the symplectic realization of Poisson-Nijenhuis manifolds}

\vspace{5mm}
\author{Fani Petalidou
\\
\\
\emph{Department of Mathematics}
\\
\emph{Aristotle University of Thessaloniki} \\
\emph{54124 Thessaloniki, Greece} \\
\\
\noindent\emph{E-mail: petalido@math.auth.gr}}
\date{}
\maketitle

\vskip 1 cm


\begin{tabular}{rr}
 & \hspace{48mm} \emph{Dedicated to Professor Charles-Michel Marle,}\\
  & \hspace{48mm} \emph{with my deepest admiration and respect,}\\
  &  \hspace{48mm} \emph{on the occasion of his 80\,$^{th}$ birthday.}
\end{tabular}

\vskip 15 mm

\begin{abstract}
We consider the problem of the symplectic realization of a Poisson-Nijenhuis manifold. By applying a new technique developed by M. Crainic and I. M$\check{\mathrm{a}}$rcu$\c{t}$ for the study of the above problem in the case of a Poisson manifold, we establish the existence, under a condition, of a nondegenerate Poisson-Nijenhuis structure on an open neighborhood of the zero-section of the cotangent bundle of the manifold, which symplectizes the initial structure. Additionally, we present some examples.
\end{abstract}

\vspace{5mm} \noindent {\bf{Keywords: }}{Bi-hamiltonian manifold, Poisson-Nijenhuis manifold, symplectic realization, contravariant connection, Poisson-spray.}

\vspace{3mm} \noindent {\bf{MSC (2010):}} 53D17, 53D05, 53D25, 37K10.

\section{Introduction}
A \emph{bi-Hamiltonian manifold} is a smooth manifold $M$ endowed with a pair $(\Pi_0,\Pi_1)$ of compatible Poisson structures in the sense that $\Pi_0 + \Pi_1$ is still a Poisson structure. The last condition happens if and only if the Schouten bracket of $\Pi_0$ with $\Pi_1$ vanishes. \emph{Poisson-Nijenhuis} manifolds is a particular class of bi-Hamiltonian manifolds which are characterized by the property that the pair $(\Pi_0, \Pi_1)$ possesses a Nijenhuis operator $N$ as a recursion operator. The first notion is due to F. Magri \cite{magri-1978} and the second has been introduced by F. Magri and C. Morosi \cite{magri-mo} in order to study the complete integrability of Hamiltonian dynamical systems. Franco Magri first discovered that if a dynamical system $X$ on $M$ can be written in Hamiltonian form in two different compatible ways, namely, there exist a bi-Hamiltonian structure $(\Pi_0, \Pi_1)$ on $M$ and $f_0,f_1 \in \C$ such that
\begin{equation*}
X = \Pi_0^\#(df_1) = \Pi_1^\#(df_0),
\end{equation*}
then it possesses an infinity of first integrals. More precisely, if $\Pi_0$ and $\Pi_1$ have Casimirs functions, then, they are first integrals of $X$ as well as the Casimirs of the Poisson pencil $\Pi_\lambda = \Pi_0 + \lambda \Pi_1$, $\lambda \in \R$. While, if $\Pi_0$ is nondegenerate, then the pair $(\Pi_0, \Pi_1)$ has a recursion operator $N = \Pi_1^\# \circ \Pi_0^{\#^{-1}}$ and the functions $I_k = \mathrm{Trace}(N^k)$, $k=1,\ldots, n=\displaystyle{\frac{1}{2}\dim M}$, are first integrals of $X$. Consequently, if enough of the obtained first integrals are functionally independent, then the system is completely integrable in the sense of Arnold-Liouville. This result is the origin of an increased interest in the study of bi-Hamiltonian and specially Poisson-Nijenhuis manifolds during the last 35 years. Many mathematicians have examined a plethora of problems related with these structures. Indicatively, we cite the works of Y. Kosmann-Schwarzbach and F. Magri \cite{yks-magri} where Poisson-Nijenhuis and related structures are studied under algebraic circumstances, \cite{yks-LMP} in which Y. Kosmann-Schwarzbach suggested a relation between Lie bialgebroids and Poisson-Nijenhuis structures, \cite{vai-N} and \cite{vai-Red-PN} in which I. Vaisman developed the theory of Poisson-Nijenhuis manifolds using Lie algebroids and studied the reduction problem of these manifolds, respectively. Moreover, we cite the works of I. M. Gel'fand and I. Zakharevich \cite{gel-zak}, P. J. Olver \cite{olver}, F. J. Turiel \cite{tu} and of the author \cite{petal-thesis}, \cite{petal-coimbra} where the local classification of bi-Hamiltonian and Poisson-Nijenhuis structures is considered.

One of the most important problems in Poisson geometry from the point of view of integration and quantization theory of Poisson manifolds is that of \emph{symplectic realization of a Poisson manifold} $(M,\Pi)$. It consists of constructing a surjective submersion $\Phi : (\tilde{M}, \tilde{\Pi}) \to (M, \Pi)$ from a symplectic-Poisson manifold $(\tilde{M}, \tilde{\Pi})$, i.e. $\tilde{\Pi}$ is nondegenerate, to $(M,\Pi)$ such that $\Phi$ is a Poisson map. The existence of a local symplectic realization for a given $(M, \Pi)$, i.e. in the neighborhood of a singular point of $\Pi$, and its universality was proven by A. Weinstein in \cite{wei} while the existence of a global symplectic realization for any Poisson manifold is established by A. Weinstein and his collaborators A. Coste and P. Dazord in \cite{cdw}. The same global result was obtained independently by M. Karasev \cite{karasev}.

The analogous problem in the framework of bi-Hamiltonian manifolds is expressed as follows: \emph{For a given bi-Hamiltonian manifold $(M, \Pi_0, \Pi_1)$ whose all the Poisson structures of the associated Poisson pencil $\Pi_\lambda = \Pi_0 + \lambda \Pi_1$, $\lambda \in \R$, are degenerate, construct a surjective submersion $\Phi : (\tilde{M}, \tilde{\Pi}_0, \tilde{\Pi}_1) \to (M, \Pi_0, \Pi_1)$ from a nondegenerate bi-Hamiltonian manifold $(\tilde{M}, \tilde{\Pi}_0, \tilde{\Pi}_1)$, i.e., at least one of the structures of $(\tilde{\Pi}_0, \tilde{\Pi}_1)$ is nondegenerate, to $(M, \Pi_0, \Pi_1)$ such that $\Phi$ is a Poisson map for the both pairs $(\tilde{\Pi}_0, \Pi_0)$ and $(\tilde{\Pi}_1, \Pi_1)$.} It is a difficult problem and it has been studied by the author in \cite{petal-sympl}. Her results are local and concern some special cases.

Recently, M. Crainic and I. M$\check{\mathrm{a}}$rcu$\c{t}$ have presented a new proof of the existence of a symplectic realization of a Poisson manifold $(M, \Pi)$ based on the theory of contravariant connections and of Poisson sprays. It consisted in the construction of a symplectic form on an open neighborhood of the zero-section of the cotangent bundle of $M$ \cite{cr-marcut}. By studying this paper the natural question which arises is: \emph{Can we use the M. Crainic and I. M$\check{\mathrm{a}}$rcu$\c{t}$'s new technique in order to investigate the symplectic realization problem of a degenerate bi-Hamiltonian manifold $(M, \Pi_0, \Pi_1)$?}

The purpose of this paper is to study the above question. Because of the difficulties that we have encountered in the consideration of the general case, we restricted our study in the case of Poisson-Nijenhuis manifolds and we proved the following result: \emph{For any Poisson-Nijenhuis manifold $(M, \Pi_0, N)$ endowed with a symmetric covariant connection $\nabla$ compatible with $N$, in a sense specified below, there exists a nondegenerate bi-Hamiltonian structure on an open neighborhood of the zero-section of the cotangent bundle of $M$ that symplectizes $(\Pi_0, N)$.}

The proof of the main result is given in Section \ref{section theorem} and for its presentation we follow the notation of M. Crainic and I. M$\check{\mathrm{a}}$rcu$\c{t}$'s paper. Section \ref{preliminaries} is devoted to the recall of some preliminary notions while in section \ref{section examples} we give some examples.

\section{Preliminaries}\label{preliminaries}
We first fix our notation and recall some important notions and results needed in the following. Let $M$ be a smooth $n$-dimensional manifold, we denote by $TM$ and $T^\ast M$ its tangent and cotangent bundle, respectively, by $\Gamma(TM)$ and $\Gamma(T^\ast M)$ the corresponding spaces of smooth sections of $TM$ and $T^\ast M$, and by $\C$ the space of smooth functions on $M$. The canonical projection of $T^\ast M$ onto the base $M$ is denoted by $\pi: T^\ast M \to M$. Finally, for any local coordinate system $(x^1,\ldots, x^n)$ of $M$ we denote by $(x^1,\ldots,x^n,y_1,\ldots,y_n)$ or $(x,y)$ the adapted coordinate system on $T^\ast M$.

\subsection{Poisson-Nijenhuis manifolds}
We recall some basic definitions concerning Poisson structures and we give the formal definition of a Poisson-Nijenhuis manifold.

A \emph{Poisson manifold} $(M, \Pi)$ is a smooth manifold $M$ equipped with a smooth bivector field $\Pi$ such that $[\Pi, \Pi]= 0$, where $[\cdot, \cdot]$ denotes the Schouten bracket, the unique natural extension of the Lie bracket between vector fields to multivector fields, \cite{lib-marle}, \cite{vai-b}, \cite{duf-zung}. The bivector field $\Pi$ defines a natural vector bundle morphism $\Pi^\# : T^\ast M \to TM$ whose the induced morphism on the space of smooth sections, also denoted by $\Pi^\#$, is defined, for any $\alpha, \beta \in \Gamma(T^\ast M)$, by
\begin{equation*}\label{def-diese}
\langle \beta, \Pi^\#(\alpha)\rangle = \Pi(\alpha, \beta).
\end{equation*}
The map $\Pi^\#$ is a Lie algebra homomorphism from the Lie algebra $(\Gamma(T^\ast M), [\cdot,\cdot]_\Pi)$ to the Lie algebra $(\Gamma(TM), [\cdot,\cdot])$, where the Lie bracket $[\cdot,\cdot]_\Pi$ on $\Gamma(T^\ast M)$ is given by
\begin{equation}\label{bracket-Lie-forms}
[\alpha, \beta]_\Pi = \mathcal{L}_{\Pi^\#(\alpha)}\beta - \mathcal{L}_{\Pi^\#(\beta)}\alpha - d(\Pi(\alpha,\beta)).
\end{equation}
In the particular case where $\alpha = df$, the vector field $\Pi^\#(df)$ is called the Hamiltonian vector field of $f$ with respect to $\Pi$ and it is denoted by $X_f$.

A differentiable map between two Poisson manifolds $\Phi : (M_1, \Pi_1) \to (M_2, \Pi_2)$ is called \emph{Poisson map} or \emph{Poisson morphism} if the vector bundle morphisms $\Pi_1^\# : T^\ast M_1 \to TM_1$ and $\Pi_2^\# : T^\ast M_2 \to TM_2$ satisfy, for all $x\in M_1$,
\begin{equation*}\label{Poisson-map}
\Pi_{2_{\Phi(x)}}^\# = \Phi_{\ast_x} \circ \Pi_{1_x}^\# \circ \Phi_x^\ast.
\end{equation*}

\vspace{2mm}

A \emph{Nijenhuis structure} on a manifold $M$ is a tensor field $N$ of type $(1,1)$, viewed also as a vector bundle map $N : TM \to TM$, with Nijenhuis torsion $T(N) : TM\times TM \to TM$ identically zero on $M$. This means that, for any pair $(X,Y)$ of vector fields on $M$,
\begin{equation*}
T(N)(X, Y) = [NX, NY] - N[NX, Y] - N[X, NY] + N^2[X,Y]\equiv 0.
\end{equation*}

\vspace{2mm}

A \emph{Poisson-Nijenhuis} manifold is a Poisson manifold $(M,\Pi_0)$ equipped with a compatible Nijenhuis structure $N$ in the sense that
\begin{equation}\label{cond-N-Pi0}
N\circ \Pi_0^\# = \Pi_0^\#\circ \,^tN
\end{equation}
and the Magri-Morosi's concomitant $C(\Pi_0, N)$ of $\Pi_0$ and $N$ is identically zero on $M$. The concomitant $C(\Pi_0,N)$ is a $T^\ast M$-valued bivector field on $M$ defined, for any $(\alpha, \beta)\in \Gamma(T^\ast M) \times \Gamma(T^\ast M)$, by
\begin{equation}\label{def-concomitant}
C(\Pi_0, N)(\alpha, \beta) = (\mathcal{L}_{\Pi_0^\#(\alpha)}\,^tN)\beta - (\mathcal{L}_{\Pi_0^\#(\beta)}\,^tN)\alpha + \,^tN d(\Pi_0(\alpha, \beta)) - d(\Pi_1(\alpha,\beta)).
\end{equation}
The condition (\ref{cond-N-Pi0}) and the vanishing of $C(\Pi_0,N)$ ensures that $\Pi_1$, defined by $\Pi_1^\# = N\circ \Pi_0^\#$, is a bivector field which satisfies the relation $[\Pi_0, \Pi_1] =0$. Then, the vanishing of $T(N)$ implies that $\Pi_1$ is also Poisson, \cite{joana-marle}. Thus, $(\Pi_0,\Pi_1)$ is a bihamiltonian structure.

As is well-known \cite{yks-magri}, a Poisson-Nijenhuis structure $(\Pi_0, N)$ defines a whole hierarchy $(\Pi_k)_{k\in \N}$ of Poisson structures, $\Pi_k^\# = N^k\circ \Pi_0^\#$, which are pairwise compatible, i.e., for all $k,l \in \N$, $[\Pi_k, \Pi_l]=0$. In the particular case where $N$ is nondegenerate, the hierarchy is defined for any $k\in \Z$, also.

\subsection{Lifts to the cotangent bundle}
The theory of \emph{lifts} of tensor fields from an arbitrary manifold $M$ to its cotangent bundle $T^\ast M$ is dealt with in the book by K. Yano and S. Ishihara \cite{Yano-Ishi}. In this subsection we shall present some results concerning the lifts on $T^\ast M$ of a Nijenhuis tensor and of a vector field on $M$.

\vspace{3mm}
Let $N$ be a Nijenhuis tensor field on $M$ which in a local coordinate system $(x^1, \ldots, x^n)$ of $M$ is written as $N = \nu^i_j \displaystyle{\frac{\partial}{\partial x^i}}\otimes dx^j$. We shall present two different ways of lifting $N$ on $T^\ast M$. In the first way, the \emph{vertical lift of} $N$, we get a vertical vector field $N^v$ on $T^\ast M$ given by $N^v = y_i \nu^i_j \displaystyle{\frac{\partial}{\partial y_j}}$. In the second way, the \emph{complete lift of} $N$, we obtain a tensor field $N^c$ of type $(1,1)$ again on $T^\ast M$; its local expression in the coordinates $(x^1, \ldots, x^n, y_1, \ldots, y_n)$ of $T^\ast M$ is
\begin{equation*}
N^c = \nu^i_j (\frac{\partial}{\partial x^i}\otimes dx^j + \frac{\partial}{\partial y_j}\otimes dy_i) + y_i(\frac{\partial \nu^i_j}{\partial x^k} - \frac{\partial \nu^i_k}{\partial x^j} )\frac{\partial}{\partial y_j}\otimes dx^k.
\end{equation*}
Moreover, $N^c$ viewed as a vector bundle map $N^c : T(T^\ast M) \to T(T^\ast M)$ has the matrix expression
\begin{equation}\label{expression N}
N^c = \begin{pmatrix} N & 0 \cr
                      A & \,^t N
\end{pmatrix},
\end{equation}
where $A = (a^j_k)$ with $a^j_k = y_i\displaystyle{(\frac{\partial\nu^i_j}{\partial x^k} - \frac{\partial\nu^i_k}{\partial x^j})}$. We have that $T(N^c) = (T(N))^c$, where $(T(N))^c$ is the complete lift on $T^\ast M$ of the skew-symmetric $(1,2)$-tensor field $T(N)$ on $M$ \cite{Yano-Ishi}. Thus we conclude
\begin{proposition}
The complete lift $N^c$ of a Nijenhuis operator $N$ on $M$ is a Nijenhuis operator on $T^\ast M$ and reciprocally.
\end{proposition}

\vspace{2mm}

Now, we assume that $M$ is endowed with a classical symmetric linear connection $\nabla$. The symmetry condition of $\nabla$ means that its torsion $T_{\nabla}$ is identically zero, i.e., for any $X,Y \in \Gamma(TM)$,
\begin{equation}\label{torsion-connection}
T_{\nabla}(X,Y) = \nabla_X Y - \nabla_Y X - [X, Y] \equiv 0.
\end{equation}
In a coordinate system $(x^1,\ldots,x^n)$, the symmetry of $\nabla$ is expressed by the fact that $\Gamma_{ij}^k = \Gamma^k_{ji}$, where $\Gamma^k_{ij}$, $i,j,k = 1,\ldots,n$, are the coefficients (Christoffel symbols) of $\nabla_{\frac{\partial}{\partial x^i}}\displaystyle{\frac{\partial}{\partial x^j}} = \Gamma_{ij}^k\displaystyle{\frac{\partial}{\partial x^k}}$. For a vector field $X = \chi^i \displaystyle{\frac{\partial}{\partial x^i}}$ on $M$, we set
\begin{equation*}
X^h = \chi^i \frac{\partial}{\partial x^i} + y_k\Gamma^k_{ij}\chi^j \frac{\partial}{\partial y_i}
\end{equation*}
and we call the vector field $X^h$ the \emph{horizontal lift of} $X$ on $T^\ast M$. The horizontal lifts of all vector fields on $M$ define the \emph{horizontal bundle} $\mathcal{H}$ on $T^\ast M$ with respect to $\nabla$ which is a distribution on $T^\ast M$ complementary to its vertical distribution $\ker \pi_\ast$, namely
\begin{equation*}
T(T^\ast M) = \mathcal{H} \oplus \ker \pi_\ast,
\end{equation*}
\cite{sau}. The sections of $\mathcal{H}$ are called \emph{horizontal vector fields} of $T^\ast M$. It is evident that a horizontal lift is indeed horizontal and, conversely, a \emph{projectable horizontal vector field}, i.e., a horizontal vector field whose coefficients of $\displaystyle{\frac{\partial}{\partial x^i}}$ are functions pulled back from $M$, is the horizontal lift of its projection. $\mathcal{H}$ is involutive if $\nabla$ is of zero curvature.

\vspace{2mm}

We have that the complete lift $N^c$ of a Nijenhuis operator $N$ maps the projectable horizontal vector fields $X^h$, $X\in \Gamma(TM)$, to
\begin{equation}\label{im hor by N compl}
N^c X^h = (NX)^h + [\nabla N]_X^v,
\end{equation}
where $[\nabla N]_X^v$ is the vertical lift of the $(1,1)$-tensor field $[\nabla N]_X$ on $M$ defined, for any $Y \in \Gamma(TM)$, by
\begin{equation*}
[\nabla N]_X Y = (\nabla_X N)Y - (\nabla_Y N)X,
\end{equation*}
\cite{Yano-Ishi}. Therefore, we get that $N^c$ conserves the horizontal distribution $\mathcal{H}$ if and only if $[\nabla N]_X = 0$, for any $X\in \Gamma(TM)$, or, equivalently, for any $X, Y \in \Gamma(TM)$,
\begin{equation}\label{cond-nabla-N}
(\nabla_X N)Y - (\nabla_Y N)X = 0.
\end{equation}
In the following, equation (\ref{cond-nabla-N}) will be referred as \emph{compatibility condition between $N$ and $\nabla$}. An easy computation yields that, in local coordinates, (\ref{cond-nabla-N}) is equivalent to the system of equations
\begin{equation}\label{cond-nabla-N locally}
\frac{\partial \nu^i_j}{\partial x^k} - \frac{\partial \nu^i_k}{\partial x^j} = \Gamma^i_{jl}\nu^l_k - \Gamma^i_{kl}\nu^l_j.
\end{equation}

\subsection{Contravariant connections }
The fundamental concept of \emph{contravariant connection} in Poisson geometry appeared firstly in R. L. Fernandes' paper \cite{fer-conn} and from then it has an essential contribution in the study of the global properties of Poisson manifolds, \cite{fer-Lie Holo}, \cite{cr-fer-int-Pois}. Its definition is inspired by that of a classical covariant connection and it is based on the general philosophical principle in Poisson geometry that \emph{the cotangent bundle plays the role of the tangent bundle and the two are related by the bundle map} $\Pi^\#$.

\vspace{2mm}

Let $(E, \tau, M)$ be a vector bundle over a Poisson manifold $(M, \Pi)$ and $\Gamma(E)$ the space of the smooth sections of $E$. A \emph{contravariant connection} on $E$ is a bilinear map
\begin{eqnarray*}
\nabla : \Gamma(T^\ast M) \times \Gamma(E) & \to & \Gamma(E) \nonumber \\
(\alpha, s) & \mapsto & \nabla_\alpha s
\end{eqnarray*}
satisfying, for any $(\alpha,  s) \in \Gamma(T^\ast M) \times \Gamma(E)$ and $f\in \C$, the following properties:
\begin{equation*}
\nabla_{f \alpha} s = f \nabla_\alpha s \quad \quad \mathrm{and} \quad \quad \nabla_\alpha (fs) = f \nabla_\alpha s + (\mathcal{L}_{\Pi^\#(\alpha)}f) s.
\end{equation*}

In order to introduce the corresponding notion of \emph{contravariant derivative along a path} for a contravariant connection $\nabla$ on $E$, we define a suitable notion of \emph{cotangent path} (or $T^\ast M$-path in the sense of \cite{cr-fer-int-Lie}) as follows.

A \emph{cotangent path $a$ with base path $\gamma$} is a curve $a : [0,1] \to T^\ast M$ sitting above some curve $\gamma : [0,1] \to M$, that means that $\pi(a(t)) = \gamma(t)$, such that
\begin{equation*}
\frac{d\gamma}{dt}(t) = \Pi^\#(a(t)).
\end{equation*}

Hence, given a $T^\ast M$-path $(a, \gamma)$ and a path $u : [0,1] \to E$ on $E$ above $\gamma$, i.e., $\tau (u(t)) = \gamma(t)$, the \emph{contravariant derivative of $u$ along $a$}, denoted by $\nabla_a u$, is well defined. We choose a time-dependent section $s$ of $E$ such that $s(t, \gamma(t)) = u(t)$ and we set
\begin{equation*}
\nabla_a u (t) = \nabla_a s_t(\gamma(t)) + \frac{d s_t}{dt}(\gamma(t)).
\end{equation*}

\vspace{2mm}

In what follows we are interested in contravariant connections on the cotangent bundle $(T^\ast M, \pi, M)$ of a Poisson manifold $(M, \Pi)$ induced by a classical covariant connection $\nabla$ on $M$. We know that any such connection $\nabla$ on $(M, \Pi)$ induces a contravariant connection $\bar{\nabla}$ on $TM$ and a contravariant connection $\tilde{\nabla}$ on $T^\ast M$ which are defined, for any $\alpha, \beta \in \Gamma(T^\ast M)$ and $X\in \Gamma(TM)$, respectively, by
\begin{equation}\label{def-contravariant connection}
\bar{\nabla}_\alpha X = \Pi^\#(\nabla_X \alpha) + [\Pi^\#(\alpha), \, X] \quad \quad \mathrm{and} \quad \quad \tilde{\nabla}_\alpha \beta = \nabla_{\Pi^\#(\beta)}\alpha +[\alpha, \beta]_{\Pi},
\end{equation}
where $[\cdot, \cdot]_\Pi$ is the Lie bracket (\ref{bracket-Lie-forms}) on $\Gamma(T^\ast M)$ defined by $\Pi$. Because
\begin{equation*}
\Pi^\# : (\Gamma(T^\ast M), [\cdot, \cdot]_\Pi) \to (\Gamma(TM), [\cdot,\cdot])
\end{equation*}
is a Lie algebra homomorphism, the two connections are related by the formula
\begin{equation}
\bar{\nabla}_\alpha (\Pi^\#(\beta)) = \Pi^\#(\tilde{\nabla}_\alpha \beta).
\end{equation}
Moreover, in the case where $\nabla$ is without torsion, they are, also, related as it is indicated in the following Lemma.
\begin{lemma}(\cite{cr-marcut})
Let $\nabla$ be a linear torsion-free connection on $M$ and $a : [0,1] \to T^\ast M$ a cotangent path with base path $\gamma$. Then, for any smooth paths $\theta : [0,1] \to T^\ast M$ and $u: [0,1] \to TM$, both above $\gamma$, the following identity holds:
\begin{equation}\label{lemma 1.3 cr-mrc}
\langle \tilde{\nabla}_{a_t} \theta_t, \, u_t\rangle + \langle \theta_t, \, \bar{\nabla}_{a_t} u_t\rangle = \frac{d}{dt}\langle \theta_t, u_t\rangle.
\end{equation}
\end{lemma}

We close this subsection by calculating the coefficients $\tilde{\Gamma}^{ij}_k$ of $\tilde{\nabla}$ in a local coordinate system $(x^1,\ldots,x^n)$ of $M$. Let $\Gamma_{ij}^k$ be the Christoffel symbols of $\nabla$ in these coordinates. We consider the extension of $\nabla$ on $T^\ast M$, also denoted by $\nabla$, and we calculate its coefficients in $(x^1,\ldots,x^n)$:
\begin{equation}\label{conn-ext-T*M}
\nabla_{\frac{\partial}{\partial x^j}}dx^i = -\Gamma_{jk}^i dx^k.
\end{equation}
Thus,
\begin{eqnarray*}
\tilde{\nabla}_{dx^i} dx^j & \stackrel{(\ref{def-contravariant connection})}{=} & \nabla_{\Pi^\#(dx^j)}dx^i +[dx^i, dx^j]_{\Pi} \, \stackrel{(\ref{bracket-Lie-forms})}{=} \, \nabla_{\Pi^{jl}\frac{\partial}{\partial x^l}}dx^i + \frac{\partial \Pi^{ij}}{\partial x^k}dx^k \\
 & = & \Pi^{jl}\nabla_{\frac{\partial}{\partial x^l}}dx^i + \frac{\partial \Pi^{ij}}{\partial x^k}dx^k \, \stackrel{(\ref{conn-ext-T*M})}{=} \, \Pi^{jl}(-\Gamma_{lk}^idx^k)+ \frac{\partial \Pi^{ij}}{\partial x^k}dx^k \\
 & = & (\Gamma_{lk}^i \Pi^{lj}+ \frac{\partial \Pi^{ij}}{\partial x^k})dx^k.
\end{eqnarray*}
Hence we get
\begin{equation}\label{Christ - T*M}
\tilde{\Gamma}_k^{ij} =  \Gamma^i_{kl}\Pi^{lj}+ \frac{\partial \Pi^{ij}}{\partial x^k}.
\end{equation}

\subsection{Poisson sprays}\label{section poisson spray}
The contravariant analogue of the classical notion of a \emph{spray} (\cite{lang}) is the one of \emph{Poisson spray} that is defined in \cite{cr-marcut} as follows.
\begin{definition}
A \emph{Poisson spray} on a Poisson manifold $(M, \Pi)$ is a vector field $\mathcal{V}_\Pi$ on $T^\ast M$ that satisfies the following properties:
\begin{enumerate}
\item[1)]
For any $\xi \in T^\ast M$, $\xi = (x,y)$,
\begin{equation}\label{property-1}
\pi_{{\ast}_\xi} (\mathcal{V}_{\Pi_\xi}) = \Pi^\#(\xi);
\end{equation}
\item[2)]
\begin{equation}\label{property-2}
m_{{t_\ast}_\xi}(\mathcal{V}_{\Pi_\xi}) = \frac{1}{t}\mathcal{V}_{\Pi_{m_t(\xi)}},
\end{equation}
where $m : \R^\ast \times T^\ast M \to T^\ast M$ is the action by dilatation of $\R^\ast$ on the fibers of $T^\ast M$, i.e., for any $\xi = (x,y) \in T^\ast M$, $m_t(\xi) = (x, ty)$.
\end{enumerate}
\end{definition}

By studying the above properties of $\mathcal{V}_\Pi$ we conclude that, in a local coordinate system $(x,y)$ of $T^\ast M$, $\mathcal{V}_\Pi$ is written as
\begin{equation}\label{expression Poisson-spray}
\mathcal{V}_\Pi(x,y) = \sum_{i,j}\Pi^{ij}(x)y_i\frac{\partial}{\partial x^j} + \sum_{i,j,k}F_k^{ij}(x)y_i y_j\frac{\partial}{\partial y_k}, \quad \quad F_k^{ij}\in \C.
\end{equation}

In the following, we describe the local behavior of the flow $\varphi$ of a Poisson spray $\mathcal{V}_\Pi$ on $(M, \Pi)$ in the neighborhood of the points of the zero-section of $T^\ast M$ because it plays a fundamental role in the proof of our main result (Theorem \ref{THEOREM}). We note that, at each $x\in M$, the zero-section corresponds the zero-covector $0_x = (x,0)$ of $T^\ast_x M$, and, at such points we have $T_{0_x}(T^\ast M) = T_x M \oplus T_{0_x}(T_x^\ast M)\cong T_xM \oplus T_x^\ast M$, since $T_{0_x}(T_x^\ast M)$ is canonically identified with $T_x^\ast M$. Therefore, each element $u$ of $T_{0_x}(T^\ast M)$ can be written as $u = (\bar{u}, \theta_u)$, where $\bar{u} = \pi_{{\ast}_{0_x}}(u)$ and $\theta_u$ is the projection of $u$ on $T_{0_x}(T_x^\ast M) \cong T_x^\ast M$. On the other hand, taking into account the local expression (\ref{expression Poisson-spray}) of $\mathcal{V}_\Pi$, we have that the points $0_x = (x,0)$, $x\in M$, are singular points of $\mathcal{V}_{\Pi}$. So, the maximal integral curves of $\mathcal{V}_{\Pi}$ through these points are the same ones, i.e., $\varphi_t(0_x) = 0_x$, for all $t\in \R$ and all $x\in M$. Hence, $\varphi_t$ is well defined on a neighborhood of the zero-section of $T^\ast M$ for all $t\in \R$ and in particular for $t\in [0,1]$. Its tangent map at these points, $\varphi_{t_{{\ast}_{0_x}}} : T_{0_x}(T^\ast M) \to T_{0_x}(T^\ast M)$, is defined by $\varphi_{t_{{\ast}_{0_x}}} = \exp(t \dot{\mathcal{V}}_{{\Pi}_{0_x}})$. In the last relation $\dot{\mathcal{V}}_{{\Pi}_{0_x}}$ is the composition of the tangent map $\mathcal{V}_{{\Pi}_{\ast_{0_x}}} : T_{0_x}(T^\ast M) \to T_{0_{0_x}}(T(T^\ast M))$ of $\mathcal{V}_{\Pi}$ at $0_x$, as soon as it is viewed as a smooth section of $T(T^\ast M)$, with the canonical projection of $T_{0_{0_x}}(T(T^\ast M))$, which is identified with $T_{0_x}(T^\ast M) \oplus T_{0_{0_x}}(T_{0_x}(T^\ast M)) \cong T_{0_x}(T^\ast M) \oplus T_{0_x}(T^\ast M)$, on its second (vertical) summand $T_{0_x}(T^\ast M)$. (For more details, see Proposition 22.3 in \cite{abr-rob}.) From (\ref{expression Poisson-spray}), we obtain that, in the coordinates $(x,y)$, $\mathcal{V}_{{\Pi}_{\ast_{0_x}}}$ has the matrix expression
\begin{equation*}
\mathcal{V}_{{\Pi}_{\ast_{0_x}}} = \begin{pmatrix} I & 0 \cr
                                                     0 & I \cr
                                                     0 & \Pi_x \cr
                                                     0 & 0 \end{pmatrix}.
\end{equation*}
Consequently, for any $u = (\bar{u}, \theta_u) \in T_{0_x}(T^\ast M)$, $\dot{\mathcal{V}}_{{\Pi}_{0_x}}(u)=\dot{\mathcal{V}}_{{\Pi}_{0_x}}(\bar{u}, \theta_u) = (-\Pi^\# (\theta_u), 0)$ and
\begin{equation}\label{im flow u_0}
\varphi_{t_{{\ast}_{0_x}}}(u) = \exp \big(t \dot{\mathcal{V}}_{{\Pi}_{0_x}}(u)\big) = (\bar{u} - t\Pi_0^\# (\theta_u), \theta_u).
\end{equation}

\vspace{2mm}

Finally, we remark that the contravariant analogue of the classical notion of \emph{geodesic spray} can also be defined in the framework of Poisson geometry. Given a contravariant connection $\nabla$ on $T^\ast M$ of $(M, \Pi)$ with coefficients $\Gamma^{ij}_k$, the notion of \emph{geodesics of} $\nabla$ on $T^\ast M$ is defined, as usual, as the cotangent paths $a : [0,1] \to T^\ast M$ whose contravariant derivative along itself is identically zero: $\nabla_{a} a(t) = 0$. In local coordinates $(x,y)$, a geodesic can be regarded as a curve $(\gamma, a) : [0,1] \to M\times T^\ast M$, $(\gamma(t), a(t)) = (x^1(t), \ldots, x^n(t), y_1(t), \ldots, y_n(t))$, which satisfies the following system of ode's:
\begin{equation}\label{system ode's}
\left\{
\begin{array}{l}
\frac{dx^i}{dt} = \Pi^{ki}(x(t))y_k(t)\\
\\
\frac{dy_i}{dt} = - \Gamma_i^{jk}(x(t))y_j(t)y_k(t)
\end{array}, \quad \quad
i = 1,\ldots,n.
\right.
\end{equation}
The system (\ref{system ode's}) determines the vector field $\mathcal{V}_\Pi$ on $T^\ast M$, that is given by
\begin{equation*}
\mathcal{V}_\Pi = \Pi^{ki}y_k\frac{\partial}{\partial x^i} - \Gamma_i^{jk}y_jy_k \frac{\partial}{\partial y_i}
\end{equation*}
and it is called the \emph{geodesic Poisson spray}, which, clearly, has the properties (\ref{property-1}) and (\ref{property-2}) of Poisson sprays. The above discussion ensures us the existence of Poisson sprays on $(M,\Pi)$.

We have that, if the considered contravariant connection on $T^\ast M$ is the connection $\tilde{\nabla}$ defined by a linear symmetric connection $\nabla$ on $TM$ as in (\ref{def-contravariant connection}), the corresponding geodesic Poisson spray takes the form
\begin{eqnarray}\label{Geodesic-Poisson-spray}
\mathcal{V}_{\Pi} & = & \sum_{i,j}\Pi^{ij}y_i\frac{\partial}{\partial x^j} - \sum_{i,j,k}\tilde{\Gamma}_k^{ij}y_i y_j\frac{\partial}{\partial y_k} \nonumber \\
&  \stackrel{(\ref{Christ - T*M})}{=} & \sum_{i,j}\Pi^{ij}y_i\frac{\partial}{\partial x^j} - \sum_{i,j,k,l}(\Gamma^i_{kl}\Pi^{lj}+ \frac{\partial \Pi^{ij}}{\partial x^k})y_i y_j\frac{\partial}{\partial y_k} \nonumber \\
& = & \sum_{i,j}\Pi^{ij}y_i\frac{\partial}{\partial x^j} - \sum_{i,j,k,l}\Gamma^i_{kl}\Pi^{lj}y_i y_j\frac{\partial}{\partial y_k} \nonumber \\
& = & \sum_i y_i (\Pi_0^\#(dx^i))^h.
\end{eqnarray}
Note that, for any $k = 1, \ldots, n$, the sum $\sum_{i,j}\displaystyle{\frac{\partial \Pi^{ij}}{\partial x^k}y_i y_j}$ is annulled because of the skew-symmetry of $\Pi$. In this special case, at every point $\xi \in T^\ast M$, $\mathcal{V}_{\Pi_\xi}$ coincides with the horizontal lift of $\Pi^\#(\xi)$ on $T^\ast M$ with respect to $\nabla$. Therefore, $\mathcal{V}_{\Pi}$ is a section of the horizontal subbundle $\mathcal{H}$ of $T(T^\ast M)$ defined by $\nabla$.

\subsection{Symplectic realization of a Poisson manifold}\label{realization Crainic-Marcut}
In this subsection we present, briefly, the basic steps of the proof of M. Crainic and I. M$\check{\mathrm{a}}$rcu$\c{t}$' theorem \cite{cr-marcut} on the \emph{symplectic realization of a Poisson manifold:}
\begin{theorem}\label{theorem-crainic-marcut}(\cite{cr-marcut})
Given a Poisson manifold $(M, \Pi)$ and a Poisson spray $\mathcal{V}_\Pi$, there exists an open neighborhood $\mathcal{U}\subset T^\ast M$ of the zero-section so that
\begin{equation*}
\Omega: = \int_0^1 \varphi_t^\ast \omega_{can}dt,
\end{equation*}
where $\varphi$ is the flow of $\mathcal{V}_\Pi$, is a symplectic structure on $\mathcal{U}$ and the canonical projection $\pi : (\mathcal{U}, \Omega) \to (M, \Pi)$ is a symplectic realization.
\end{theorem}

Firstly, they calculate the value of $\Omega$ on vectors tangent to $T^\ast M$ at zeros $0_x \in T_x^\ast M$. By identifying the space $T_{0_x}(T^\ast M)$ with $T_x M \oplus T_x^\ast M$, they prove that, for all $u=(\bar{u}, \theta_u), w=(\bar{w}, \theta_w) \in T_{0_x}(T^\ast M) \cong T_x M \oplus T_x^\ast M$,
\begin{equation}\label{values Omega 0}
\Omega_{0_x}(u,w) = \langle \theta_w, \bar{u}\rangle - \langle \theta_u, \bar{w}\rangle - \Pi (\theta_u, \theta_w).
\end{equation}
The above formula implies that the closed $2$-form $\Omega$ is nondegenerate at the points of the zero-section of $T^\ast M$. So, there exists a neighborhood $\mathcal{U}$ of the zero-section in $T^\ast M$ on which $\varphi_t$ is well defined for any $t\in [0,1]$ and $\Omega \vert_{\mathcal{U}}$ is symplectic. In the second step, by fixing such a $\mathcal{U}$, they consider a torsion-free covariant connection $\nabla$ on $TM$, in order to handle vectors tangent to $T^\ast M$, and they establish a generalization of (\ref{values Omega 0}) at arbitrary points $\xi \in \mathcal{U}$. Precisely, they prove (Lemma 2.3 in \cite{cr-marcut}) that, for any pair $(u,w)$ of elements of $T_{\xi}(T^\ast M)$,
\begin{equation*}
\Omega(u,w) = \big(\langle \tilde{\theta}_{w_t}, \bar{u}_t\rangle   - \langle \tilde{\theta}_{u_t}, \bar{w}_t \rangle - \Pi (\tilde{\theta}_{u_t}, \tilde{\theta}_{w_t}) \big)\big\vert_0^1,
\end{equation*}
where $u_t = \varphi_{t_{\ast}}(u)$ (resp. $w_t = \varphi_{t_{\ast}}(w)$), $\bar{u}_t = \pi_\ast(u_t)$ (resp. $\bar{w}_t = \pi_\ast(w_t)$), $\theta_{u_t} = u_t - \bar{u}_t^h$ (resp. $\theta_{w_t} = w_t - \bar{w}_t^h$) and $\tilde{\theta}_{u_t}$ (resp. $\tilde{\theta}_{w_t}$) is a solution of the differential equation $\tilde{\nabla}_{a_t}\tilde{\theta}_{u_t} = \theta_{u_t}$ (resp. $\tilde{\nabla}_{a_t}\tilde{\theta}_{w_t} = \theta_{w_t}$), with $a_t$ being the path on $\mathcal{U}$ given by $a_t = \varphi_t(\xi)$. In the third and last step, they show that the projection $\pi\vert_{\mathcal{U}} : \mathcal{U} \to M$ push-down the bivector associated to $\Omega$ to the Poisson tensor $\Pi$, i.e., that $\pi$ is a Poisson map.

\section{Sympectic realization of Poisson-Nijenhuis manifolds}\label{section theorem}
Let $(M, \Pi_0, N)$ be a Poisson-Nijenhuis manifold. Without loss of generality we assume that $N$ is nondegenerate. We can make this assumption because, if $\det N = 0$, we can replace $N$ with the nondegenerate Nijenhuis operator $N'= I + N$ which produces with $\Pi_0$ the same bi-Hamiltonian structure viewed as a bi-parametric family $\Pi_{\kappa, \lambda} =\kappa \Pi_0 + \lambda \Pi_1$, $\kappa, \lambda \in \R$, of pairwise compatible Poisson structures, with $\Pi_1^\# = N\circ \Pi_0^\#$. In this case, an hierarchy $(\Pi_k)_{k \in \Z}$, $\Pi_k^\# = N^k\circ \Pi_0^\#$, of pairwise compatible Poisson structures is also defined on $M$. We consider the complete lift $N^c$ of $N$ on $T^\ast M$ and the pair of sympectic forms $(\omega_{can}, \omega_1)$, where $\omega_{can}$ is the canonical symplectic form on $T^\ast M$ and $\omega_1$ is the symplectic $2$-form defined by $\omega_1 (\cdot, \cdot) = \omega_{can} (N^c \cdot, \cdot) = \omega_{can}( \cdot, N^c \cdot)$. In the local coordinate system $(x,y)$ of $T^\ast M$, they have, respectively, the matrix expression
\begin{equation*}
\omega_{can} = \begin{pmatrix} 0 & -I \cr
                      I & 0
\end{pmatrix}
\quad \quad \mathrm{and} \quad \quad \omega_1 = \begin{pmatrix} -A & - \,^tN \cr
                      N & 0
\end{pmatrix}.
\end{equation*}
Since $N^c$ is a Nijenhuis operator, $\omega_{can}$ and $\omega_1$ are Poisson-compatible in the sense of \cite{tu, magri-mo}. Furthermore, we consider a Poisson spray $\mathcal{V}_{\Pi_0}$ on $T^\ast M$ associated to $\Pi_0$ and we denote by $\varphi$ its flow. Thereafter, we endow $T^\ast M$ with the pair of closed $2$-forms
\begin{equation}\label{Omega0-Omega1}
\Omega_0 = \int_0^1 \varphi_t^\ast \omega_{can} dt \quad \quad \mathrm{and} \quad \quad \Omega_1 = \int_0^1 \varphi_t^\ast \omega_1 dt
\end{equation}
and we remember that $\Omega_0$ is symplectic \cite{cr-marcut}.

\begin{lemma}\label{lemma compatible}
The pair $(\Omega_0,\Omega_1)$ is a pair of Poisson-compatible $2$-forms.
\end{lemma}
\begin{proof}
By definition \cite{tu}, $(\Omega_0, \Omega_1)$ is Poisson-compatible if its recursion operator $R$ defined by $R = \Omega_0^{\flat^{-1}}\circ \Omega_1^\flat$ is a Nijenhuis operator.\footnote{We recall that $\Omega_i^\flat$, $i=0,1$, denotes the vector bundle map from $TM$ to $T^\ast M$ whose the induced map on the space of smooth sections, also denoted by $\Omega_i^\flat$, is defined as follows: for any $X,Y \in \Gamma(TM)$, $\langle \Omega_i^\flat(X),Y\rangle = - \Omega_i(X,Y)$.} It is well known that it is true if and only if the 2-form $\Omega_2$ defined by the relation $\Omega_2^\flat = \Omega_1^\flat \circ R = \Omega_0^\flat \circ R^2$ is closed, \cite{br}. We have that, for all $u,w \in \Gamma(T(T^\ast M))$,
\begin{equation*}
\Omega_1(u,w) = \int_0^1 (\varphi_t^\ast \omega_1)(u,w) dt = \int_0^1  \omega_1(\varphi_{t \ast}u, \, \varphi_{t \ast}w) dt = \int_0^1  \omega_{can}(N^c (\varphi_{t \ast}u), \, \varphi_{t \ast}w)
\end{equation*}
and
\begin{equation*}
\Omega_1(u,w) =  \Omega_0(Ru, w) = \int_0^1 \varphi_t^\ast \omega_{can} (Ru, w) dt = \int_0^1 \omega_{can} ( \varphi_{t_\ast}Ru, \varphi_{t_\ast}w) dt.
\end{equation*}
Thus, for all $u,w \in \Gamma(T(T^\ast M))$,
\begin{equation}\label{two expressions of Omega-2}
\int_0^1 \omega_{can} ( \varphi_{t_\ast}Ru, \varphi_{t_\ast}w) dt = \int_0^1  \omega_{can}(N^c (\varphi_{t_\ast}u), \, \varphi_{t_\ast}w).
\end{equation}
We consider on $T^\ast M$ the closed $2$-form $\omega_2$ defined by $\omega_2^\flat = \omega_0^\flat \circ (N^c)^2$ (since $T(N^c)=0$, $d\omega_2 =0$) and we calculate, for all $u,w \in \Gamma(T(T^\ast M))$,
\begin{eqnarray*}
\Omega_2(u,w) & = & \Omega_0(R^2u,w) =  \int_0^1 \varphi_t^\ast \omega_{can} (R^2u, w) dt = \int_0^1 \omega_{can} ( \varphi_{t_\ast}R^2u, \varphi_{t_\ast}w) dt \nonumber \\
              & \stackrel{(\ref{two expressions of Omega-2})}{=} & \int_0^1  \omega_{can}(N^c (\varphi_{t_\ast}(Ru)), \, \varphi_{t_\ast}w)dt = \int_0^1  \omega_{can}(\varphi_{t_\ast}(Ru), \, N^c \varphi_{t_\ast}w)dt \nonumber \\
              & = & \int_0^1  \omega_{can}(\varphi_{t_\ast}(Ru), \, \varphi_{t_\ast}(\varphi_{-t_\ast}N^c \varphi_{t_\ast}w))dt \nonumber \\
              & \stackrel{(\ref{two expressions of Omega-2})}{=} & \int_0^1  \omega_{can}(N^c (\varphi_{t_\ast}u), \, \varphi_{t_\ast}(\varphi_{-t_\ast}N^c \varphi_{t_\ast}w))dt \nonumber \\
              & = & \int_0^1  \omega_{can}(N^c (\varphi_{t_\ast}u), \, N^c (\varphi_{t_\ast}w))dt = \int_0^1  \omega_{can}((N^c)^2 (\varphi_{t_\ast}u), \, \varphi_{t_\ast}w)dt \nonumber \\
              & = & \int_0^1  \omega_{2}(\varphi_{t_\ast}u, \, \varphi_{t_\ast}w)dt = \int_0^1 \varphi_t^\ast \omega_2(u,w) dt.
\end{eqnarray*}
Therefore,
\begin{equation*}
\Omega_2 = \int_0^1 \varphi_t^\ast \omega_2 dt \quad \quad \mathrm{and} \quad \quad d\Omega_2 =0.
\end{equation*}
Hence we get the compatibility of $\Omega_0$ with $\Omega_1$.
\end{proof}

\vspace{2mm}

For the follow of our study we have also need the next lemma.

\begin{lemma}\label{lemma nabla N}
Let $(M,\Pi_0,N)$ be a Poisson-Nijenhuis manifold equipped with a torsion-free covariant connection $\nabla$ compatible with $N$. Then, for any $\alpha, \beta \in \Gamma(T^\ast M)$, the following identity holds:
\begin{equation}\label{cond-tilde-nabla-N}
\tilde{\nabla}_{\alpha}(\,^tN \beta) = \,^tN(\tilde{\nabla}_{\alpha}\beta),
\end{equation}
where $\tilde{\nabla}$ is the contravariant connection (\ref{def-contravariant connection}) on $T^\ast M$ induced by $\nabla$ and $\Pi_0$.
\end{lemma}
\begin{proof}
Effectively, for any $\alpha, \beta \in \Gamma(T^\ast M)$ and $X\in \Gamma(TM)$, we have
\begin{eqnarray*}
\langle \tilde{\nabla}_{\alpha}(\,^tN \beta), X\rangle & \stackrel{(\ref{def-contravariant connection})}{=} & \langle \nabla_{\Pi_0^\#(\,^tN \beta)}\alpha + [\alpha, \,^tN\beta]_{\Pi_0},\, X  \rangle \nonumber \\
& \stackrel{(\ref{bracket-Lie-forms})}{=} & \langle \nabla_{\Pi_0^\#(\,^tN \beta)}\alpha , X\rangle + \langle \mathcal{L}_{\Pi_0^\#(\alpha)}(\,^tN \beta) - \mathcal{L}_{\Pi_0^\#(\,^tN \beta)}\alpha -d(\Pi_1(\alpha,\beta)),\, X\rangle \nonumber \\
&\stackrel{(\ref{cond-N-Pi0})}{=} & \Pi_1^\#(\beta) \langle \alpha, \, X\rangle - \langle \alpha, \nabla_{N\Pi_0^\#(\beta)}X\rangle + \Pi_0^\#(\alpha)\langle \,^tN \beta,\, X\rangle \nonumber \\
& &  -\, \langle \,^tN \beta,\, \mathcal{L}_{\Pi_0^\#(\alpha)}X\rangle -\Pi_1^\#(\beta)\langle \alpha, \,X\rangle + \langle \alpha, \, \mathcal{L}_{\Pi_1^\#(\beta)}X\rangle \nonumber \\
& & - \,\langle d(\Pi_1(\alpha,\beta)),\, X\rangle \nonumber \\
& \stackrel{(\ref{torsion-connection})}{=} & -\, \langle \alpha, [N\Pi_0^\#(\beta), X] + \nabla_X(N\Pi_0^\#(\beta))\rangle + \Pi_0^\#(\alpha)\langle \,^tN \beta,\, X\rangle \nonumber \\
& & -\, \langle \,^tN \beta,\, \mathcal{L}_{\Pi_0^\#(\alpha)}X\rangle +  \langle \alpha, [\Pi_1^\#(\beta), X]\rangle - \langle d(\Pi_1(\alpha,\beta)),\, X\rangle \nonumber \\
& = & - \,\langle \alpha, (\nabla_XN)\Pi_0^\#(\beta) + N \nabla_X \Pi_0^\#(\beta)\rangle + \Pi_0^\#(\alpha)\langle \,^tN \beta,\, X\rangle  \nonumber \\
& & - \, \langle \,^tN \beta,\, \mathcal{L}_{\Pi_0^\#(\alpha)}X\rangle - \langle d(\Pi_1(\alpha,\beta)),\, X\rangle.
\end{eqnarray*}
On the other hand,
\begin{eqnarray*}
\langle \,^tN (\tilde{\nabla}_{\alpha}\beta), \, X\rangle & \stackrel{(\ref{def-contravariant connection})}{=} & \langle \nabla_{\Pi_0^\#(\beta)}\alpha + [\alpha, \beta]_{\Pi_0}, \, NX\rangle \nonumber \\
& \stackrel{(\ref{bracket-Lie-forms})}{=} & \langle \nabla_{\Pi_0^\#(\beta)}\alpha + \mathcal{L}_{\Pi_0^\#(\alpha)}\beta - \mathcal{L}_{\Pi_0^\#(\beta)}\alpha - d(\Pi_0(\alpha, \beta)), \, NX  \rangle \nonumber \\
& = & \Pi_0^\#(\beta)\langle \alpha, \,NX\rangle - \langle \alpha, \,\nabla_{\Pi_0^\#(\beta)}(NX)\rangle + \Pi_0^\#(\alpha)\langle \beta, \, NX\rangle  \nonumber \\
&  &-\, \langle \beta, \, \mathcal{L}_{\Pi_0^\#(\alpha)}(NX)\rangle - \Pi_0^\#(\beta) \langle \alpha, NX\rangle + \langle \alpha, \mathcal{L}_{\Pi_0^\#(\beta)}(NX)\rangle  \nonumber \\
& & -\, \langle d(\Pi_0(\alpha, \beta)),\, NX\rangle \nonumber \\
& = & -\, \langle \alpha, \, (\nabla_{\Pi_0^\#(\beta)}N)X + N \nabla_{\Pi_0^\#(\beta)}X\rangle + \Pi_0^\#(\alpha)\langle \beta, \, NX\rangle \nonumber \\
& & -\, \langle \beta, \, (\mathcal{L}_{\Pi_0^\#(\alpha)}N)X + N \mathcal{L}_{\Pi_0^\#(\alpha)}X \rangle + \langle \alpha, \mathcal{L}_{\Pi_0^\#(\beta)}N(X) + N \mathcal{L}_{\Pi_0^\#(\beta)}X\rangle  \nonumber \\
& & - \, \langle d(\Pi_0(\alpha, \beta)),\, NX\rangle \nonumber \\
& = & - \langle \alpha, \, (\nabla_{\Pi_0^\#(\beta)}N)X \rangle - \langle \,^tN\alpha, \, \nabla_{\Pi_0^\#(\beta)}X\rangle  + \Pi_0^\#(\alpha)\langle \beta, \, NX\rangle \nonumber \\
& &-  \,\langle (\mathcal{L}_{\Pi_0^\#(\alpha)}\,^tN)\beta, \, X \rangle  - \langle \,^tN \beta, \, \mathcal{L}_{\Pi_0^\#(\alpha)}X\rangle + \langle (\mathcal{L}_{\Pi_0^\#(\beta)}\,^tN)\alpha, \,X\rangle \nonumber \\
& & + \,\langle \,^tN \alpha, \, \mathcal{L}_{\Pi_0^\#(\beta)}X\rangle -\langle \,^tN d(\Pi_0(\alpha, \beta)),\, X\rangle.
\end{eqnarray*}
Hence, taking into account the facts that $N$ is compatible with the symmetric connection $\nabla$ and $C(\Pi_0,N)=0$, we obtain
\begin{eqnarray*}
\langle \tilde{\nabla}_{\alpha}(\,^tN \beta) - \,^tN (\tilde{\nabla}_{\alpha}\beta), \,X \rangle & = & - \langle \alpha, (\nabla_XN)\Pi_0^\#(\beta) - (\nabla_{\Pi_0^\#(\beta)}N)X \rangle \nonumber \\
& & + \langle \,^tN \alpha, \nabla_{\Pi_0^\#(\beta)}X - \nabla_X \Pi_0^\#(\beta) - [\Pi_0^\#(\beta), \,X]\rangle \nonumber \\
& & + \,\langle - d(\Pi_1(\alpha,\beta)) + (\mathcal{L}_{\Pi_0^\#(\alpha)}\,^tN)\beta - (\mathcal{L}_{\Pi_0^\#(\beta)}\,^tN)\alpha \nonumber \\
& & + \,^tN d(\Pi_0(\alpha, \beta)),\, X\rangle \nonumber \\
&\stackrel{(\ref{cond-nabla-N}, \ref{torsion-connection}, \ref{def-concomitant})}{=}& 0,
\end{eqnarray*}
for any $X\in \Gamma(TM)$. So, equality (\ref{cond-tilde-nabla-N}) is established. The above result means that, under the conditions of compatibility of $N$ with $(\Pi_0, \nabla)$, $\,^t N$ is parallel with respect to $\tilde{\nabla}$.
\end{proof}

\vspace{2mm}
Now, we proceed with the proof of the central theorem of this work.

\begin{theorem}\label{THEOREM}
Let $(M, \Pi_0, N)$ be a Poisson-Nijenhuis manifold, with $N$ nondegenerate, equipped with a torsion-free covariant connection $\nabla$ compatible with $N$. Let, also, $(\Pi_k)_{k\in \Z}$, $\Pi_k^\# = N^k\circ \Pi_0^\#$, be the associated hierarchy of pairwise compatible Poisson structures on $M$, $\mathcal{V}_{\Pi_0}$ a Poisson spray corresponding to $\Pi_0$ and $\varphi$ its flow. Then, there exists an open neighborhood $\mathcal{U}$ of the zero-section in $T^\ast M$ such that the canonical projection $\pi : (\mathcal{U}, \Omega_0, \Omega_1) \to (M, \Pi_0, \Pi_{-1})$ is a symplectic realization of $(M, \Pi_0, \Pi_{-1})$, where $(\Omega_0, \Omega_1)$ is the pair of Poisson-compatible symplectic structures on $T^\ast M$ defined by (\ref{Omega0-Omega1}).
\end{theorem}

\begin{proof}
The result of M. Crainic and I. M$\check{\mathrm{a}}$rcu$\c{t}$ \cite{cr-marcut} (see, also, subsection \ref{realization Crainic-Marcut}) ensures the existence of a neighborhood $\mathcal{U}$ in $T^\ast M$ of the zero-section of $T^\ast M$ such that $\varphi_t$ is defined for all $t\in [0,1]$, $\Omega_0\vert_{\mathcal{U}}$ is symplectic and $\pi : (\mathcal{U}, \tilde{\Pi}_0) \to (M,\Pi_0)$, where $\tilde{\Pi}_0$ is the symplectic Poisson structure defined by $\Omega_0$, is a Poisson map. Our aim is to prove the claims that $\Omega_1\vert_{\mathcal{U}}$ is symplectic and $\pi : (\mathcal{U}, \tilde{\Pi}_1) \to (M,\Pi_{-1})$ is also a Poisson map, where $\tilde{\Pi}_1$ is the symplectic Poisson structure defined by $\Omega_1$.

\vspace{2mm}
\noindent
\emph{First step:} We start by evaluating $\Omega_1$ on vectors tangent to $T^\ast M$ at the points $0_x = (x,0)$ of the zero-section of $T^\ast M$. As we noted in the subsection \ref{section poisson spray}, $T_{0_x}(T^\ast M) = T_x M \oplus T_{0_x}(T_x^\ast M)$ is canonically identified with $T_xM \oplus T_x^\ast M$ and each element $u$ of $T_{0_x}(T^\ast M)$ is identified with $u = (\bar{u}, \theta_u)$, where $\bar{u} = \pi_{{\ast}_{0_x}}(u)$ and $\theta_u$ is the projection of $u$ on $T_{0_x}(T_x^\ast M) \cong T_x^\ast M$. Hence, for any pair $(u,w)$ of elements of $T_{0_x}(T^\ast M)$, $u = (\bar{u}, \theta_u)$ and $w = (\bar{w}, \theta_w)$, we have
\begin{eqnarray}\label{omega1_0}
\omega_{1_{0_x}}(u,w)& = &\begin{pmatrix}\bar{w} & \theta_w \end{pmatrix}\begin{pmatrix} A & -\,^tN \cr
                      N & 0 \end{pmatrix}_{0_x}\begin{pmatrix} \bar{u} \cr \theta_u \end{pmatrix} = \begin{pmatrix}\bar{w} & \theta_w \end{pmatrix}\begin{pmatrix} 0 & -\,^tN \cr
                      N & 0 \end{pmatrix}\begin{pmatrix} \bar{u} \cr \theta_u \end{pmatrix} \nonumber \\
                      & = & \langle \theta_w, N\bar{u}\rangle - \langle \,^tN\theta_u, \bar{w}\rangle = \langle \theta_w, N\bar{u}\rangle - \langle \theta_u, N\bar{w}\rangle.
\end{eqnarray}
Furthermore, since $\varphi_t(0_x) = 0_x$,
\begin{eqnarray*}
(\varphi_t^\ast \omega_1)_{0_x}(u,w) & = & \omega_{1_{0_x}} \big((\varphi_t)_{{\ast}_{0_x}}(u), (\varphi_t)_{{\ast}_{0_x}}(w) \big)\nonumber \\
& \stackrel{(\ref{im flow u_0})}{=}&  \omega_{1_{0_x}}\big((\bar{u} - t\Pi_0^\# (\theta_u), \theta_u),(\bar{w} - t\Pi_0^\# (\theta_w), \theta_w) \big) \nonumber \\
& \stackrel{(\ref{omega1_0})}{=} & \langle \theta_w, N\bar{u} - tN\Pi_0^\# (\theta_u)\rangle - \langle \theta_u, N\bar{w}  - tN\Pi_0^\# (\theta_w)\rangle \nonumber \\
& = &  \langle \theta_w, N\bar{u}\rangle   - \langle \theta_u, N\bar{w} \rangle - 2t\Pi_1 (\theta_u, \theta_w).
\end{eqnarray*}
Consequently,
\begin{eqnarray}\label{Omega_1 - 0}
\Omega_{1_{\,0_x}}(u,w)& = & \int_0^1 (\varphi_t^\ast \omega_1)_{0_x}(u,w)dt = \int_0^1 \big(\langle \theta_w, \bar{u}\rangle   - \langle \theta_u, \bar{w} \rangle - 2t\Pi_1 (\theta_u, \theta_w) \big)dt  \nonumber \\
& = & \langle \theta_w, N\bar{u}\rangle   - \langle \theta_u, N\bar{w} \rangle - \Pi_1 (\theta_u, \theta_w).
\end{eqnarray}
The last expression implies that $\Omega_1$ is nondegenerate at the points $0_x$, $x\in M$, of $T^\ast M$. Since $\varphi_t$ is defined on $\mathcal{U}$ for any $t\in [0,1]$, we conclude that $\Omega_1$ is symplectic on $\mathcal{U}$.

\vspace{2mm}
\noindent
\emph{Second step:} In this step we evaluate $\Omega_1$ on vectors tangent to $T^\ast M$ at arbitrary points $\xi$ of $\mathcal{U}$ and we establish a formula analogous of (\ref{Omega_1 - 0}) which we shall use in the proof of the assertion that $\pi : (\mathcal{U}, \Omega_1) \to (M,\Pi_{-1})$ is a Poisson morphism. In order to describe the sections of $T(T^\ast M)$, we assume that $M$ is equipped with a symmetric covariant connection $\nabla$ compatible with $N$ in the sense of (\ref{cond-nabla-N}).\footnote{For some comments on the existence of a such connection, see Remark {\ref{remark-existence-connection}}.} Then, the tangent bundle of $T^\ast M$ is decomposed, with respect to $\nabla$, as $T(T^\ast M) = \mathcal{H}\oplus \ker \pi_{\ast}$, where $\mathcal{H}$ is the horizontal distribution on $T^\ast M$ defined by $\nabla$. Hence, any tangent vector $u$ of $T^\ast M$ at $\xi$ is written as
\begin{equation*}
u = \bar{u}^h + \theta_u,
\end{equation*}
where $\bar{u}^h$ is the horizontal lift at $\xi$ of the projection $\bar{u} = \pi_{\ast_{\xi}}(u)$ of $u$ on $T_{\pi(\xi)}M$ and $\theta_u$ is the projection of $u$ on $\ker \pi_{{\ast}_\xi} \cong T^\ast_{\pi(\xi)}M$ parallel to $\bar{u}^h$. Clearly, at the points $\xi = 0_x$, the decomposition $T_\xi(T^\ast M) = \mathcal{H}_\xi\oplus \ker \pi_{{\ast}_\xi}$ coincides with the one described in the previous step. The fact that $\nabla$ is torsion-free ensures that the distribution $\mathcal{H}$ is Lagrangian with respect to $\omega_{can}$. The extra condition (\ref{cond-nabla-N}) of compatibility of $N$ with $\nabla$ implies that $\mathcal{H}$ is conserved by $N^c$. Hence, $\mathcal{H}$ is also Lagrangian with respect to $\omega_1$, i.e., $\mathcal{H}$ is a bi-Lagrangian distribution with respect to $(\omega_{can}, \omega_1)$. Indeed, for all $u,w \in T_\xi(T^\ast M)$,
\begin{equation}\label{ypolo Nc}
N^cu = N^c\bar{u}^h + N^c\theta_u  \stackrel{(\ref{im hor by N compl},\ref{cond-nabla-N})}{=} (N\bar{u})^h + \,^tN\theta_u,
\end{equation}
and
\begin{eqnarray}\label{omega1_xi}
\omega_1(u,w) & = & \omega_{can}(N^cu, w) \nonumber \\
              & \stackrel{(\ref{ypolo Nc})}{=} & \omega_{can}((N\bar{u})^h + \,^tN\theta_u , \bar{w}^h + \theta_w) \nonumber \\
              & = & \omega_{can}((N\bar{u})^h, \bar{w}^h) + \omega_{can}((N\bar{u})^h, \theta_w) + \omega_{can}(\,^tN\theta_u, \bar{w}^h) + \omega_{can}(\,^tN\theta_u, \theta_w) \nonumber \\
              & = & \langle \theta_w, N\bar{u}\rangle - \langle \,^tN\theta_u, \bar{w}\rangle = \langle \theta_w, N\bar{u}\rangle - \langle \theta_u, N\bar{w}\rangle,
\end{eqnarray}
because $\mathcal{H}$ and $\ker \pi_\ast$ are Lagrangian distributions on $T^\ast M$ with respect to $\omega_{can}$. The above formula is a generalization of (\ref{omega1_0}) at an arbitrary point $\xi \in T^\ast M$.

In the following, we fix $\xi$ in $\mathcal{U}$ and we consider the cotangent path $a : [0,1] \to \mathcal{U}$, $a_t : = \varphi_t(\xi)$, which is the integral curve of $\mathcal{V}_{\Pi_0}$ through $\xi$, and we denote by $\gamma = \pi \circ a$ its base path on $M$. By pushing forward a tangent vector $u \in T_\xi \mathcal{U}$ by $\varphi_{{t_\ast}_\xi} : T_\xi\mathcal{U} \to T_{a_t}\mathcal{U}$, $t\in [0,1]$, we obtain a smooth path $u_t : = \varphi_{{t_\ast}_\xi}(u)$ of vectors along $a$; its projection $\bar{u}_t = \pi_{{\ast}_{a_t}}(u_t)$ on $TM$ yields a path of vectors along $\gamma$ while its projection $\theta_{u_t}$ on the vertical space $\ker \pi_{{\ast}_{a_t}}$ parallel to $\bar{u}_t^h$ defines a path of covectors along $\gamma$, also. I.e., we have $\bar{u}_t \in T_{\gamma(t)}M$ and $\theta_{u_t} \in T^\ast_{\gamma(t)}(M)$. Accordingly to Lemma 2.2 of \cite{cr-marcut}, the two paths are related by
\begin{equation}\label{lemma 2.2 cr-mrc}
\bar{\nabla}_{a_t}\bar{u}_t = \Pi_0^\#(\theta_{u_t}).
\end{equation}
The action of $N^c$ on $u_t$ produces another path of vectors along $a$ given by
\begin{equation*}
N^c u_t \stackrel{(\ref{ypolo Nc})}{=} (N\bar{u}_t)^h + \,^tN\theta_{u_t}.
\end{equation*}
Its corresponding paths on $TM$ and $T^\ast M$ are $N\bar{u}_t$ and $\,^tN\theta_{u_t}$, respectively. Hence, taking into account (\ref{lemma 2.2 cr-mrc}), we get
\begin{equation}\label{rel bar nabla Pi1}
\bar{\nabla}_{a_t}(N \bar{u}_t) = \Pi_0^\#(\,^tN\theta_{u_t}) = \Pi_1^\#(\theta_{u_t}).
\end{equation}

Now, we can establish a generalization of (\ref{Omega_1 - 0}) at an arbitrary point $\xi$ of $\mathcal{U}$. We consider a pair $(u,w)$ of elements of $T_\xi \mathcal{U}$ and its associated pairs of paths $(u_t, w_t)$ on $T\mathcal{U}$ over $a$, $(\bar{u}_t, \bar{w}_t)$ on $TM$ and $(\theta_{u_t}, \theta_{w_t})$ on $T^\ast M$, both over the base path $\gamma$ of $a$. Let $\tilde{\theta}_{u_t}$ (resp. $\tilde{\theta}_{w_t}$) be a path in $T^\ast M$ solution of the differential equation
\begin{equation}\label{diff-equation theta}
\tilde{\nabla}_{a_t}\tilde{\theta}_{u_t} = \theta_{u_t} \quad \quad (\mathrm{resp.} \quad \tilde{\nabla}_{a_t}\tilde{\theta}_{w_t} = \theta_{w_t}).
\end{equation}
We will show that
\begin{equation}\label{expression Omega 1 xi}
\Omega_1(u,w) = \big(\langle \tilde{\theta}_{w_t}, N\bar{u}_t\rangle   - \langle \tilde{\theta}_{u_t}, N\bar{w}_t \rangle - \Pi_1 (\tilde{\theta}_{u_t}, \tilde{\theta}_{w_t}) \big)\big\vert_0^1.
\end{equation}
We have
\begin{equation*}
\Omega_1(u,w) = \int_0^1 (\varphi_t^\ast \omega_1)(u,v)dt = \int_0^1 \omega_1 (u_t, w_t)dt \stackrel{(\ref{omega1_xi})}{=} \int_0^1 \big(\langle \theta_{w_t}, N\bar{u}_t\rangle - \langle \theta_{u_t}, N\bar{w}_t\rangle \big)dt.
\end{equation*}
Hence, it is enough to prove that
\begin{equation}\label{1}
\langle \theta_{w_t}, N\bar{u}_t\rangle - \langle \theta_{u_t}, N\bar{w}_t\rangle = \frac{d}{dt}\big( \langle \tilde{\theta}_{w_t}, N\bar{u}_t\rangle   - \langle \tilde{\theta}_{u_t}, N\bar{w}_t \rangle - \Pi_1 (\tilde{\theta}_{u_t}, \tilde{\theta}_{w_t}) \big).
\end{equation}
In fact, we have
\begin{eqnarray*}
\langle \theta_{w_t}, \, N\bar{u}_t\rangle - \langle \theta_{u_t}, \, N\bar{w}_t\rangle & = & \langle \tilde{\nabla}_{a_t}\tilde{\theta}_{w_t}, \, N\bar{u}_t  \rangle - \langle \tilde{\nabla}_{a_t}\tilde{\theta}_{u_t}, \, N\bar{w}_t  \rangle  \nonumber \\
& \stackrel{(\ref{lemma 1.3 cr-mrc})}{=} & \frac{d}{dt} \big(\langle \tilde{\theta}_{w_t}, \, N\bar{u}_t \rangle - \langle \tilde{\theta}_{u_t}, \, N\bar{w}_t\rangle \big) \nonumber \\
& & - \langle \tilde{\theta}_{w_t}, \,\bar{\nabla}_{a_t}(N \bar{u}_t)  \rangle + \langle \tilde{\theta}_{u_t}, \,\bar{\nabla}_{a_t}(N \bar{w}_t)   \rangle.
\end{eqnarray*}
Taking into account (\ref{rel bar nabla Pi1}) and Lemma \ref{lemma nabla N}, the last two terms yield
\begin{eqnarray*}
\langle \tilde{\theta}_{w_t}, \,\bar{\nabla}_{a_t}(N \bar{u}_t)  \rangle - \langle \tilde{\theta}_{u_t}, \,\bar{\nabla}_{a_t}(N \bar{w}_t)   \rangle & = &  \langle \tilde{\theta}_{w_t}, \, \Pi_1^\#(\theta_{u_t})\rangle - \langle \tilde{\theta}_{u_t}, \, \Pi_1^\#(\theta_{w_t})\rangle \nonumber \\
& = &  \langle \tilde{\theta}_{w_t}, \, N\Pi_0^\#(\tilde{\nabla}_{a_t}\tilde{\theta}_{u_t})\rangle - \langle \tilde{\theta}_{u_t}, \, \Pi_0^\# \,^tN(\tilde{\nabla}_{a_t}\tilde{\theta}_{w_t} )\rangle \nonumber \\
& = & \langle \,^tN \tilde{\theta}_{w_t}, \,\bar{\nabla}_{a_t}\Pi_0^\#(\tilde{\theta}_{u_t})\rangle + \langle \,^tN(\tilde{\nabla}_{a_t}\tilde{\theta}_{w_t} ), \, \Pi_0^\#(\tilde{\theta}_{u_t})\rangle \nonumber \\
& \stackrel{(\ref{cond-tilde-nabla-N})}{=} &  \langle \,^tN \tilde{\theta}_{w_t}, \,\bar{\nabla}_{a_t}\Pi_0^\#(\tilde{\theta}_{u_t})\rangle + \langle  \tilde{\nabla}_{a_t}(\,^tN\tilde{\theta}_{w_t} ), \, \Pi_0^\#(\tilde{\theta}_{u_t})\rangle \nonumber \\
& = &  \frac{d}{dt} \langle \,^tN \tilde{\theta}_{w_t}, \, \Pi_0^\#(\tilde{\theta}_{u_t}) \rangle \nonumber \\
& = &  \frac{d}{dt}(\Pi_1 (\tilde{\theta}_{u_t}, \, \tilde{\theta}_{w_t})).
\end{eqnarray*}
Thus, relation (\ref{1}) is true and, consequently, formula (\ref{expression Omega 1 xi}) is also satisfied.

\vspace{2mm}
\noindent
\emph{Third step:} In the third and last step, we shall prove that the projection $\pi$ is a Poisson map for the pair $(\tilde{\Pi}_1, \Pi_{-1})$. For this, we firstly remark that
\begin{equation*}
\mathrm{orth}_{\Omega_1}(\ker \pi_\ast) = \mathrm{orth}_{\Omega_0}(R\ker \pi_\ast) \quad \quad \mathrm{and} \quad \quad R\mathrm{orth}_{\Omega_1}(\ker \pi_\ast) = \mathrm{orth}_{\Omega_0}(\ker \pi_\ast),
\end{equation*}
where $\mathrm{orth}_{\Omega_i}(\cdot)$ is the orthogonal distribution of $(\cdot)$ with respect to the symplectic structure $\Omega_i$, $i=0,1$, and $R$ is the recursion operator of $(\Omega_0,\Omega_1)$. We will show that
\begin{equation}\label{equation orth}
\mathrm{orth}_{\Omega_1}(\ker \pi_\ast) = \ker(\pi_1)_\ast,
\end{equation}
where $\pi_1 = \pi\circ \varphi_1$. Since $\varphi_1$ is a diffeomorphism of $\mathcal{U}$, $\dim \ker(\pi_1)_\ast = \dim \ker \pi_\ast =n$ and $\dim \mathrm{orth}_{\Omega_1}(\ker \pi_\ast) =n$. So, it suffices to show that $\ker(\pi_1)_\ast \subseteq \mathrm{orth}_{\Omega_1}(\ker \pi_\ast)$. We fix $\xi \in \mathcal{U}$ and we consider a vector $u\in \ker \pi_{\ast_\xi}$ and a vector $w \in \ker(\pi_1)_{\ast_\xi}$, then $\bar{u}_0 = \pi_{\ast_\xi} (u)=(\pi \circ \varphi_0)_{\ast_\xi} (u)=0$, since $\varphi_0 = id$, and $\bar{w}_1 =0$. In view of $N^c \ker \pi_\ast \subseteq \ker \pi_\ast$, we have $(\overline{N^cu})_0 = N\bar{u}_0=0$. On the other hand, we remark that the differential equation (\ref{diff-equation theta}), as an equation on $\tilde{\theta}_{u_t}$ (resp. $\tilde{\theta}_{w_t}$), is a linear ordinary differential equation having solutions defined for any $t\in [0,1]$ and satisfying any given initial condition. So, we can choose solutions satisfying the conditions $\tilde{\theta}_{u_0}=0$ and $\tilde{\theta}_{w_1} =0$. Thus,
\begin{eqnarray*}
\Omega_1(u,w) & \stackrel{(\ref{expression Omega 1 xi})}{=} & \big( \langle \tilde{\theta}_{w_t}, N\bar{u}_t\rangle   - \langle \tilde{\theta}_{u_t}, N\bar{w}_t \rangle - \Pi_1 (\tilde{\theta}_{u_t}, \tilde{\theta}_{w_t}) \big)\vert_0^1 \nonumber \\
              & = & \langle \tilde{\theta}_{w_1}, N\bar{u}_1\rangle   - \langle \tilde{\theta}_{u_1}, N\bar{w}_1 \rangle - \Pi_1 (\tilde{\theta}_{u_1}, \tilde{\theta}_{w_1}) \nonumber \\
              &   & -\langle \tilde{\theta}_{w_0}, N\bar{u}_0\rangle   + \langle \tilde{\theta}_{u_0}, N\bar{w}_0 \rangle + \Pi_1 (\tilde{\theta}_{u_0}, \tilde{\theta}_{w_0}) \nonumber \\
              & = & 0,
\end{eqnarray*}
whence we conclude that $w \in \mathrm{orth}_{\Omega_{1_\xi}}(\ker \pi_{\ast_\xi})$ and that (\ref{equation orth}) is valid. But, $\ker(\pi_1)_\ast = \mathrm{orth}_{\Omega_0}(\ker \pi_\ast)$ (\cite{cr-marcut}) and $\mathrm{orth}_{\Omega_1}(\ker \pi_\ast) = \mathrm{orth}_{\Omega_0}(R\ker \pi_\ast)$, therefore, $\mathrm{orth}_{\Omega_0}(R\ker \pi_\ast) = \mathrm{orth}_{\Omega_0}(\ker \pi_\ast)$, which means that $\ker \pi_\ast$ is invariant by $R$.

Now we will calculate the projection of $\tilde{\Pi}_1$ by $\pi$ using the two expressions of $\Omega_1$ (see, Lemma \ref{lemma compatible}). We consider a point $\xi\in \mathcal{U}$ and an arbitrary covector $\theta \in T_x^\ast M$, where $x=\pi(\xi)$, and we denote by $u$ the unique vector in $T_\xi \mathcal{U}$ defined by the relation
\begin{equation*}
\Omega_1^\flat(u) = \pi^\ast \theta.
\end{equation*}
We can easily remark that $u$ is a point of $\mathrm{orth}_{\Omega_{1_\xi}}(\ker \pi_\ast{_\xi}) = \ker(\pi_1)_{\ast_\xi}$, thus $\bar{u}_1 =0$. Furthermore, $\Omega_1^\flat(u) = \pi^\ast \theta \Leftrightarrow \Omega_0^\flat(Ru) = \pi^\ast \theta$, which yields that $Ru \in \mathrm{orth}_{\Omega_{0_\xi}}(\ker \pi_{\ast_\xi}) = \ker(\pi_1)_{\ast_\xi}$, so $(\overline{Ru})_1 =0$. Consequently, for any $w\in T_\xi \mathcal{U}$, we have
\begin{eqnarray*}
\Omega_1(u,w) & = & \Omega_0(Ru,w) = \big( \langle \tilde{\theta}_{w_t}, (\overline{Ru})_t\rangle   - \langle \tilde{\theta}_{(Ru)_t}, \bar{w}_t \rangle - \Pi_0 (\tilde{\theta}_{(Ru)_t}, \tilde{\theta}_{w_t}) \big)\vert_0^1 \nonumber \\
              & = & \langle \tilde{\theta}_{w_1}, (\overline{Ru})_1\rangle   - \langle \tilde{\theta}_{(Ru)_1}, \bar{w}_1 \rangle - \Pi_0 (\tilde{\theta}_{(Ru)_1}, \tilde{\theta}_{w_1}) \nonumber \\
              &   & -\langle \tilde{\theta}_{w_0}, (\overline{Ru})_0\rangle + \langle \tilde{\theta}_{(Ru)_0}, \bar{w}_0 \rangle + \Pi_0 (\tilde{\theta}_{(Ru)_0}, \tilde{\theta}_{w_0})
\end{eqnarray*}
and, by choosing $\tilde{\theta}_{(Ru)_t}$ verifying the initial condition $\tilde{\theta}_{(Ru)_1} =0$, we obtain
\begin{eqnarray*}
\Omega_1(u,w) & = & -\langle \tilde{\theta}_{w_0}, (\overline{Ru})_0\rangle   + \langle \tilde{\theta}_{(Ru)_0}, \bar{w}_0 \rangle + \Pi_0 (\tilde{\theta}_{(Ru)_0}, \tilde{\theta}_{w_0})\, \Leftrightarrow \nonumber \\
\langle \theta, \bar{w}_0\rangle & = & -\langle \tilde{\theta}_{w_0}, (\overline{Ru})_0\rangle   + \langle \tilde{\theta}_{(Ru)_0}, \bar{w}_0 \rangle + \Pi_0 (\tilde{\theta}_{(Ru)_0}, \tilde{\theta}_{w_0}).
\end{eqnarray*}
The last equation holds for any $w\in T_\xi \mathcal{U}$ and for any initial value of $\tilde{\theta}_w$. Thus, $\tilde{\theta}_{(Ru)_0} = \theta$ and $\Pi_0^\#(\theta) = (\overline{Ru})_0$. On the other hand, from the second expression of $\Omega_1$ we take
\begin{eqnarray*}
\Omega_1(u,w) & = & \big( \langle \tilde{\theta}_{w_t}, N\bar{u}_t\rangle   - \langle \tilde{\theta}_{u_t}, N\bar{w}_t \rangle - \Pi_1 (\tilde{\theta}_{u_t}, \tilde{\theta}_{w_t}) \big)\vert_0^1 \nonumber \\
              & = & \langle \tilde{\theta}_{w_1}, N\bar{u}_1\rangle   - \langle \tilde{\theta}_{u_1}, N\bar{w}_1 \rangle - \Pi_1 (\tilde{\theta}_{u_1}, \tilde{\theta}_{w_1}) \nonumber \\
              &   & -\langle \tilde{\theta}_{w_0}, N\bar{u}_0\rangle   + \langle \tilde{\theta}_{u_0}, N\bar{w}_0 \rangle + \Pi_1 (\tilde{\theta}_{u_0}, \tilde{\theta}_{w_0}).
\end{eqnarray*}
Then, by choosing $\tilde{\theta}_{u_1} =0$ and because of $N\bar{u}_1 =0$, we get
\begin{equation*}
\langle \theta, \bar{w}_0\rangle = -\langle \tilde{\theta}_{w_0}, N\bar{u}_0\rangle   + \langle \,^tN\tilde{\theta}_{u_0}, \bar{w}_0 \rangle + \Pi_0 (\,^tN\tilde{\theta}_{u_0}, \tilde{\theta}_{w_0}),
\end{equation*}
which holds for any $w\in T_\xi \mathcal{U}$ and for any initial value of $\tilde{\theta}_w$. Therefore, $\,^tN\tilde{\theta}_{u_0} = \theta$ and $\Pi_0^\#(\,^tN\tilde{\theta}_{u_0}) = N\bar{u}_0$. Hence,
\begin{equation}\label{Pi-1}
\Pi_0^\#(\theta) = \Pi_0^\#(\,^tN\tilde{\theta}_{u_0})= N\bar{u}_0 \Leftrightarrow N^{-1}\Pi_0^\#(\theta) = \bar{u}_0 \Leftrightarrow \Pi_{-1}^\#(\theta) = \bar{u}_0.
\end{equation}
However,
\begin{eqnarray*}
\pi^\ast \theta = \Omega_1^\flat(u) & \Leftrightarrow & \tilde{\Pi}_1^\# (\pi^\ast \theta) = u \,\Rightarrow  \, \pi_\ast (\tilde{\Pi}_1^\# (\pi^\ast \theta)) = \pi_\ast(u) \nonumber \\
& \Leftrightarrow & (\pi_\ast \circ \tilde{\Pi}_1^\# \circ \pi^\ast)(\theta)=\bar{u}_0 \, \stackrel{(\ref{Pi-1})}{\Leftrightarrow} \, (\pi_\ast \circ \tilde{\Pi}_1^\# \circ \pi^\ast)(\theta) = \Pi_{-1}^\#(\theta).
\end{eqnarray*}
Therefore,
\begin{equation*}
\pi_\ast \circ \tilde{\Pi}_1^\# \circ \pi^\ast = \Pi_{-1}^\#,
\end{equation*}
which means that $\pi$ is a Poisson map for the pair $(\tilde{\Pi}_1, \Pi_{-1})$, also.
\end{proof}

\begin{remark}
\emph{Let $\tilde{\Pi}_k$ be the Poisson structure defined, for any $k\in \Z$, by the symplectic form}
\begin{equation*}
\Omega_k = \int_0^1 \varphi_t^\ast \omega_k dt, \quad \quad \mathrm{where} \quad \quad \omega_k(\cdot,\cdot) = \omega_{can}((N^c)^k\cdot,\cdot).
\end{equation*}
\emph{By considering on $\mathcal{U}$ the hierarchy $(\tilde{\Pi}_k)_{k \in \Z}$ of pairwise compatible Poisson structures, we can easily prove that $\pi$ is a Poisson map for any pair $(\tilde{\Pi}_k, \Pi_{-k})$, $k\in \Z$.}
\end{remark}

\begin{remark}\label{remark-existence-connection}
\emph{In the second step of the proof of the above theorem, we have assumed the existence of a symmetric covariant connection $\nabla$ on $M$ compatible with $N$. In general, it is difficult to establish the conditions under which a given tensor field $S$ on a smooth manifold $M$ admits a ``compatible'', in a certain sense, symmetric covariant connection. In local coordinates, finding a torsionless $\nabla$ ``compatible'' with $S$ reduces to determining the existence of solutions for non-homogeneous $\C$-linear systems with unknowns the Christoffel symbols of $\nabla$. In our case, for given $N$, the corresponding systems are the ones given by (\ref{cond-nabla-N locally}):}
\begin{equation}\label{system gamma}
\Gamma^i_{jl}\nu^l_k - \Gamma^i_{kl}\nu^l_j = \frac{\partial \nu^i_j}{\partial x^k} - \frac{\partial \nu^i_k}{\partial x^j}.
\end{equation}
\emph{By calling, for each $i=1,\ldots,n$, $\Gamma^i$ the symmetric matrix with elements the unknown functions $\Gamma^i_{jk}$, (\ref{system gamma}) is written as}
\begin{equation*}
\Gamma^iN - \,^tN\Gamma^i = \big(\frac{\partial \nu^i_j}{\partial x^k} - \frac{\partial \nu^i_k}{\partial x^j} \big).
\end{equation*}
\emph{It is a $\C$-linear system of $\displaystyle{\frac{n^2-n}{2}}$ equations with $\displaystyle{\frac{n^2+n}{2}}$ unknowns, the functions $\Gamma^i_{jk}$. Such a system, it is either incompatible, or, if it admits a solution, then it admits an infinity of solutions whose difference is a solution of the corresponding homogeneous system.}
\end{remark}

\section{Examples}\label{section examples}
In this section we give two examples of Poisson-Nijenhuis structures for which our result is applicable.
\begin{example}
Pair of diagonal quadratic Poisson structures: \emph{Let $V$ be a finite dimensional (real) vector space and $(x^1,\ldots, x^n)$ a system of linear coordinates for $V$. We recall that any bivector field $\Pi$ on $V$ of type}
\begin{equation*}
\Pi = \sum_{i<j}\varpi^{ij}x^ix^j\frac{\partial}{\partial x^i}\wedge \frac{\partial}{\partial x^j}, \quad \quad \mathrm{with} \quad \quad \varpi^{ij} \in \R,
\end{equation*}
\emph{is Poisson and it is called \emph{diagonal quadratic Poisson structure}, \cite{duf-zung}, \cite{cam-pol-an}. Therefore, any pair $(\Pi_0, \Pi_1)$ of quadratic Poisson structures is compatible in the sense of Magri-Morosi and it defines a bi-Hamiltonian structure on $V$. We consider such a pair and we suppose that its elements are related by a recursion operator $N$, i.e., $\Pi_1^\# = N\Pi_0^\# = \Pi_0^\# \,^tN$. Precisely, if}
\begin{equation*}
\Pi_0 = \sum_{i<j}\varpi_0^{ij}x^ix^j\frac{\partial}{\partial x^i}\wedge \frac{\partial}{\partial x^j} \quad \quad \mathrm{and} \quad \quad \Pi_1 = \sum_{i<j}\varpi_1^{ij}x^ix^j\frac{\partial}{\partial x^i}\wedge \frac{\partial}{\partial x^j},
\end{equation*}
\emph{then the components $\nu^i_j$ of $N$ must be of type $\nu^i_j = n^i_j \displaystyle{\frac{x^i}{x^j}}$ with $n^i_j \in \R$ and $\varpi_1^{ij} = n^i_l \varpi_0^{lj} = \varpi_0^{il}n_l^j$. Because $(\Pi_0, \Pi_1)$ are Poisson compatible, $N$ is necessary a Nijenhuis operator. Our problem consists of finding a symmetric covariant connection $\nabla$ on $V$ such that its Christoffel symbols $\Gamma^i_{jk}$ verify condition (\ref{cond-nabla-N locally}). We can easily check that the connection with Christoffel symbols $\Gamma^i_{ii} = - \displaystyle{\frac{1}{x^i}}$ and all the other $\Gamma_{ij}^k$ zero gives a solution to our problem. Hence, we conclude that any Poisson-Nijenhuis structure of the considered type is symplectically realizable.}
\end{example}

\begin{example}
\emph{We consider the pair $(\Pi_0,\Pi_1)$ of compatible Poisson structures on $\R^6$, where $\Pi_0$ is the linear Poisson structure associated with the periodic Toda lattice of $3$-particles and $\Pi_1$ is the Poisson structure constructed in \cite{dam-pet} that has the same Casimir invariants with $\Pi_0$ (the functions $C = a_1a_2a_3$ and $C'=b_1+b_2+b_3$) and whose first part is the quadratic Poisson bracket associated to Volterra lattice. In Flaschka's coordinate system $(a_1,a_2,a_3, b_1,b_2,b_3)$,}
\begin{equation*}
\Pi_0 = \sum_{i=1}^3 a_i\frac{\partial}{\partial a_i}\wedge (\frac{\partial}{\partial b_i} - \frac{\partial}{\partial b_{i+1}}) \quad \quad \mathrm{and} \quad \quad \Pi_1 = \sum_{i=1}^3(a_ia_{i+1}\frac{\partial}{\partial a_i}\wedge \frac{\partial}{\partial a_{i+1}} + \frac{\partial}{\partial b_i}\wedge \frac{\partial}{\partial b_{i+1}}),
\end{equation*}
\emph{with the convention $(a_{i+3},b_{i+3}) = (a_i,b_i)$. The pair possesses an infinity of recursion operators of type}
\begin{equation*}
N = \begin{pmatrix} 0 & 0 & 0 & f & f & f + a_1 \cr
                    0 & 0 & 0 & g + a_2 & g & g \cr
                    0 & 0 & 0 & h & h+a_3 & h \cr
                    \displaystyle{\frac{a_3}{a_1}}F & \displaystyle{\frac{a_3}{a_2}F+\frac{1}{a_2}} & F & 0 & 0 & 0 \cr
                    \displaystyle{\frac{a_2}{a_1}G} & G & \displaystyle{\frac{a_2}{a_3}G+\frac{1}{a_3}} &  0 & 0 & 0 \cr
                    \displaystyle{\frac{a_3}{a_1}H + \frac{1}{a_1}} & \displaystyle{\frac{a_3}{a_2}H} & H & 0 & 0 & 0
\end{pmatrix},
\end{equation*}
\emph{where $f,g,h,F,G,H\in C^\infty(\R^6)$. Since, $\Pi_1^\# = N\circ \Pi_0^\#$ is Poisson and compatible with $\Pi_0$, $N$ is a Nijenhuis operator. In the case where $F=G=H=0$ and $f=g=h=C$, there exists an infinity of symmetric connections compatible with $N$. Such a $\nabla$ is defined by the functions}
\begin{eqnarray*}
\Gamma^1_{44}=\Gamma^1_{45} = \Gamma^1_{54}=\Gamma^1_{46} = \Gamma^1_{64}=\Gamma^1_{55}=\Gamma^1_{56} = \Gamma^1_{65}= C  & \mathrm{and} & \Gamma^1_{66} = C + a_1, \\
\Gamma^2_{45}= \Gamma^2_{54}=\Gamma^2_{46} = \Gamma^2_{64} =\Gamma^2_{55}=\Gamma^2_{56} = \Gamma^2_{65}=\Gamma^2_{66}= C  & \mathrm{and} & \Gamma^2_{44} = C + a_2, \\
\Gamma^3_{44}=\Gamma^3_{45} = \Gamma^3_{54} =\Gamma^3_{46}=\Gamma^3_{64} =\Gamma^3_{56} = \Gamma^3_{65} =\Gamma^3_{66}= C  & \mathrm{and} & \Gamma^3_{55} = C + a_3,
\end{eqnarray*}
\emph{and all the other $\Gamma^i_{jk}$ are zero. Then, by applying Theorem \ref{THEOREM} we conclude that $(\Pi_0, \Pi_1)$ is symplectizable.}
\end{example}

\vspace{5mm}
\noindent
\textbf{Open problem:} Our first approach to the study of the problem mentioned in the Introduction was the following. Let $(M, \Pi_0, \Pi_1)$ be a bi-Hamiltonian manifold endowed with a symmetric covariant connection $\nabla$. We consider the convex linear combination $\Pi_s = (1-s)\Pi_0 + s\Pi_1$, $s\in [0,1]$, of $\Pi_0$ and $\Pi_1$ which produces on $M$ an $1$-parameter family of pairwise compatible Poisson structures. Then the Poisson sprays and the contravariant connections corresponding to the above structures have the following nice properties.
\begin{itemize}
\item
If $\mathcal{V}_{\Pi_0}$ is a Poisson spray of $\Pi_0$ and $\mathcal{V}_{\Pi_1}$ is a Poisson spray of $\Pi_1$, then the vector field
\begin{equation*}
\mathcal{V}_{\Pi_s} = (1-s)\mathcal{V}_{\Pi_0} + s\mathcal{V}_{\Pi_1}
\end{equation*}
is a Poisson spray of $\Pi_s = (1-s)\Pi_0 + s\Pi_1$.
\item
If $\bar{\nabla}_i$ and $\tilde{\nabla}_i$, $i=0,1$, are the contravariant connections on $TM$ and $T^\ast M$, respectively, defined by the pair $(\nabla, \Pi_i)$, $i=0,1$, as in (\ref{def-contravariant connection}), then
\begin{itemize}
\item
$\bar{\nabla}_s = (1-s)\bar{\nabla}_0 + s \bar{\nabla}_1$ is a contravariant connection on $TM$,
\item
$\tilde{\nabla}_s = (1-s)\tilde{\nabla}_0 + s \tilde{\nabla}_1$ is a contravariant connection on $T^\ast M$,
\end{itemize}
and the two are related by the formula
\begin{equation*}
(\bar{\nabla}_s)_\alpha \Pi_s^\# (\beta) = \Pi_s^\#((\tilde{\nabla}_s)_\alpha\beta).
\end{equation*}
\end{itemize}
Hence, if $(\varphi_s)_t$ is the flow of $\mathcal{V}_{\Pi_s}$, by applying the M. Crainic and I. M$\check{\mathrm{a}}$rcu$\c{t}$'s technique we construct on a convenable neighborhood $\mathcal{U}$ of the $0$-section of $T^\ast M$ an $1$-parameter family of symplectic forms
\begin{equation*}
\Omega_s = \int_0^1 (\varphi_s)_t^\ast \omega_{can}dt, \quad \quad s\in [0,1],
\end{equation*}
such that $\pi : (\mathcal{U}, \tilde{\Pi}_s) \to (M, \Pi_s)$, where $\tilde{\Pi}_s = \Omega_s^{-1}$, is a Poisson map for any $s\in [0,1]$. The question which arises is: \emph{Are the Poisson structures $\tilde{\Pi}_s$, $s\in [0,1]$, compatible between them? If not, under what conditions can this happen?}

\vspace{5mm}
\noindent
\textbf{Acknowledgments}
This work is partially supported by the project DRASI A of AUTH Research Committee.

\end{document}